\documentclass[12pt,thmsa]{article}
\usepackage{amsfonts}
\usepackage{amsmath}
\usepackage{amssymb}
\usepackage{graphics}

\newtheorem{theorem}{Theorem}

\newtheorem{example}[theorem]{Example}

\newtheorem{lemma}[theorem]{Lemma}
\newtheorem{notation}[theorem]{Notation}

\newtheorem{proposition}[theorem]{Proposition}
\newtheorem{remark}[theorem]{Remark}

\newenvironment{proof}[1][Proof]{\noindent\textbf{#1.} }{\ \rule{0.5em}{0.5em}}
\input{tcilatex}
\begin{document}

\begin{center}
{\Large Finitely presented groups and homotopy of presentations of
triangular algebras.}
\end{center}

\vspace{0.3cm}

\bigskip

\begin{center}
{\ Jorge Nicol\'{a}s L\'{o}pez}

\vspace{0.3cm}{\small jnlopez@mdp.edu.ar}

\bigskip

\textit{Departamento de Matem\'{a}tica,\\[2mm]
Facultad de Ciencias Exactas y Naturales, \\[2mm]
Universidad Nacional de Mar del Plata.}
\end{center}

\bigskip

Abstract: Given any finitely presented group $G$ we find a \textbf{triangular%
} algebra such that has two presentations, one with fundamental group $G$
and another with trivial group. Thus proving that given a collection $%
G_{1},\ldots ,G_{n}$ of finitely presented groups there exist a \textbf{%
triangular} algebra $A$ such that all $G_{i}$ appear as fundamental group of
some presentation of $A$, extending one of the result of \cite%
{BustaCastonguay}.

\section{Introduction}

Let $A$ be a basic, connected, finite dimensional algebra over a closed
field $k$. There exists a unique quiver $Q$ and a two-side ideal $I$ of the
path algebra $kQ,$ such that $A=kQ/I$. Since, in general, $I$ is not unique,
the morphism $v:kQ\rightarrow A=kQ/I$, as well as $(Q,I)$ is called a 
\textbf{presentation} of $A$. Each presentation has associated a \textbf{%
fundamental group} $\pi _{1}(Q,I)$. In \cite{fiscbacher} it is proved that
given a finitely presented group $G$, there exists a bounded quiver $(Q,I)$
such that $\pi _{1}(Q,I)=G$. In \cite{BustaCastonguay} this results it is
generalized, namely: for any collection of finitely presented groups $%
\{G_{i}\}_{i=1}^{n}$ there exists an algebra $A$ and a collections of
presentations $\{(Q,I_{i})\}_{i=1}^{n}$ of this algebra such that $\pi
_{1}(Q,I_{i})=G_{i}$ for $i=1,\ldots ,n$. Also it is proved a similar result
by \textbf{triangular} algebras, but, in this case, the involved groups are
obtained from cyclic groups by performing finite free and direct products.
In the present work we leave this restriction over the groups and, thus, we
proved the Theorem:

\begin{theorem}
\label{theo unico}Let $G_{1},\ldots ,G_{n}$ be finitely presented groups.
Then there exists a \textbf{triangular} algebra $A$ having presentations $(%
\hat{Q},I_{i})$ for $i=1,\ldots ,n$, such that $\pi _{1}(\hat{Q}%
,I_{i})=G_{i} $.
\end{theorem}

To prove this we build a special quiver $\hat{Q}$ and show the appropriated
ideals $I_{i}$. In the first place, for any finitely presented group $G_{i}$
we built a quiver $Q_{G_{i}}$ and two presentations $(Q_{G_{i}},J_{i})$ and $%
(Q_{G_{i}},\bar{J}_{i})$ such that $\pi _{1}(Q_{G_{i}},J_{i})=G_{i}$ and $%
\pi _{1}(Q_{G_{i}},\bar{J}_{i})=\{1\}$. In the second place, we make a 
\textbf{co-product}, (topologically, it is a type of amalgamate of bounded
quivers) 
\begin{equation}
(\hat{Q},I_{i}):=\left( \dbigsqcup\limits_{j\neq i}(Q_{G_{j}},\bar{J}%
_{j})\right) \sqcup (Q_{G_{i}},J_{i}).  \label{for amalgamagrande}
\end{equation}%
Is not difficult prove that the fundamental group of a co-product is the
free product of the fundamental groups of each bonded quivers (see \cite%
{Castonguay}). Thus we have $\pi _{1}(\hat{Q},I_{i})=(\sqcup _{j\neq
i}\{1\})\sqcup G_{i}=G_{i}$. This is basically the same way made in \cite%
{Castonguay}. Our principal achievement is obtain each $Q_{G_{i}}$ be
triangular and find the appropriated $J_{i}$ and $\bar{J}_{i}$. The
constructions of $Q_{G_{i}}$ and $J_{i}$ are based over the presentations of
the groups by generators and relations.

In section \ref{sec fgwfg}, given a group $G$, we write $G=%
\mathbb{Z}
^{\sqcup m}\sqcup H$ where $H$ has not free generator. After that find a
convenient presentation $\left\langle S^{\prime \prime }|R^{\prime \prime
}\right\rangle $ of $H$. The relations of $R^{\prime \prime }$ are the kind
of $w_{i}=w_{j}$ or $w_{k}\ldots w_{l}=1$. In section \ref{sec quiver and
ideal} we used this presentation to built a quiver $Q_{H}$ and ideals $I$
and $\bar{I}$ such that $\pi _{1}(Q_{H},I)=H$ and $\pi _{1}(Q_{H},\bar{I}%
)=\{1\}$. The quiver $Q_{H}$ is essentially a chain of triangles $T_{i}$
where each triangle represents an element of $S^{\prime \prime }$ (Sec. \ref%
{subsec The quiver}). Each relation in $R^{\prime \prime }$ is in
correspondence with a minimal relation in $I$, in this approximated way: if $%
w_{i}$ is involved in the relation of $R^{\prime \prime }$, then $T_{i}$ is
involved in the minimal relation (\ref{subsubsec the ideal for h}). To get $%
\bar{I}$, the idea is kill all triangles. For this we built $\bar{I}$ as the
image of $I$ by a automorphism (Sec. \ref{subsubsec the ideal for 1}). In
the Sec. \ref{subsec presentation for G}, to recuperated $G$, we make
appropriates co-product with the bounder quivers of the Lemma \ref{lemmaZm},
getting the mentioned $(Q_{G_{i}},J_{i})$ and $(Q_{G_{i}},\bar{J}_{i})$.

\section{Preliminaries\label{sec Preliminares}}

\subsection{Quivers and algebras\label{subsec Quivers and Algebras}}

A \textbf{quiver} $Q$ is a quadruple $(Q_{0},Q_{1},s,t)$ where $Q_{0}$ and $%
Q_{1}$ are sets and them elements are called \textbf{vertices} and \textbf{%
arrows} respectively; $s,t$ $:Q_{1}\rightarrow Q_{0}$ are functions which
indicate the source and the target of each arrow respectively. A \textbf{path%
} $w$ is a sequence of arrows $w=\alpha _{1}\alpha _{2}\cdots \alpha _{n}$
such that $t(\alpha _{i})=s(\alpha _{i+1})$ for $i=1,\ldots n$. When $%
s(\alpha _{1})=x$ and $t(\alpha _{n})=y$, the fact is symbolized $%
w:x\rightarrow y$, and two paths with this proprieties are called \textbf{%
parallel}. A quiver $Q$ is said \textbf{finite} if $Q_{0}$ and $Q_{1}$ are
finite, and $Q$ is said \textbf{connected} when the underlying graph is
connected. We will only consider finite and connected quivers.

Given a commutative field $k$ and a quiver $Q$, the path algebra $kQ$ is the 
$k$-vector space whose base is the set of paths of $Q$, including one
stationary path $e_{x}$ for each vertex $x$ of $Q$. The multiplication of
two basis elements of $kQ$ is their composition whenever it is possible, and 
$0$ in otherwise. Let $F$ be the two-sided ideal of kQ generated by the
arrows of $Q$. A two-sided ideal $I$ of $kQ$ is called \textbf{admissible }%
if there exists an integer $m\geq 2$ such that $F^{m}\subset I\subset F^{2}$%
. The pair $(Q,I)$ is called a \textbf{bound quiver}. When a distinguished
vertex $x\in Q_{0}$ is consider, $(Q,I,x)$ is called \textbf{pointed bound
quiver.}

Conversely, if $A$ is a basic, connected and finite dimensional algebra over 
$k$, then there exists a unique finite connected quiver $Q$ and a surjective
morphism of $k$-algebras $\nu :kQ\rightarrow A$, which is not unique, with $%
I=\ker \nu $ an admissible ideal (see \cite{Gebriel} and \cite{bongGabriel}%
). The morphism $\nu $ and the pair $(Q,I)$ are called \textbf{presentations}
of the algebra $A$. Remark that a morphism $\nu $ $:kQ\rightarrow A$ is a
presentation of $A$ whenever $\{\nu (e_{x})|x\in Q_{0}\}$ is a complete set
of primitive orthogonal idempotents and, for any fixed $x,y\in Q_{0}$, we
have that $\{\nu (\alpha )+\mathrm{rad}^{2}A|\alpha :x\rightarrow y\in
Q_{1}\}$ is basis of $\nu (e_{x})(\mathrm{rad}A/\mathrm{rad}^{2}A)\nu
(e_{y}) $. An algebra $A$ is \textbf{triangular} whenever $Q$ has not
oriented cycles. In this work we will say that the quiver is triangular as
well. For further reference of bound quivers in the representation theory of
algebras can see \cite{AssemSimonSko} and \cite{aunlander}.

\begin{remark}
\label{remark Tranversiones}Let $\nu $ be a morphisms $\nu :kQ\rightarrow
kQ/I\simeq A$ defined by $\nu (e_{x})=e_{x}+I$ for $x\in Q_{0}$, and, $\nu
(\alpha )=\alpha +\rho _{\alpha }+I$ for $\alpha \in Q_{1}$ where $\rho
_{\alpha }$ is a linear combination of paths parallel to $\alpha $ and the
paths have length at least $2$, then $(Q,\ker \nu )$ is a presentation of $A$%
. In particular if one only arrow $\beta $ is transformed and $\rho _{\beta
} $ has one term, then $\nu $ is a \textbf{transvection}.
\end{remark}

\subsection{Fundamental group of a bounder quiver\label{subsec fundam group}}

Given a bound quiver $(Q,I)$, its fundamental group is defined as follows
(see \cite{Martinez}). For $x,y\in Q_{0}$, is defined the set $%
I(x,y)=e_{x}(kQ)e_{y}\cap I$. A relation $\rho =\sum_{i=1}^{m}\lambda
_{i}w_{i}\in I(x,y)$ (where $\lambda _{i}\in k^{\ast }$, and $w_{i}$ are
different paths from $x$ to $y$) is said to be \textbf{minimal} if $m\geq 2$%
, and, for every proper subset $J$ $\subset \{1,\ldots ,m\}$, we have $%
\sum_{i\in J}\lambda _{i}w_{i}\notin I(x,y)$. For a given arrow $\alpha
:x\rightarrow y$, let $\alpha ^{-1}:x\rightarrow y$ be its formal inverse. A 
\textbf{walk} $w$ in $Q$ from $x$ to $y$ is a composition $w=\alpha
_{1}^{\varepsilon _{1}}\alpha _{2}^{\varepsilon _{2}}\cdots \alpha
_{n}^{\varepsilon _{n}}$ such that the $s(\alpha _{1}^{\varepsilon _{1}})=x$%
, $t(\alpha _{n}^{\varepsilon _{n}})=y$, and, $s(\alpha _{i}^{\varepsilon
_{i}})=t(\alpha _{i-1}^{\varepsilon _{i-1}})$ for $i=1,\ldots ,n$. Define
the homotopic relation $\sim $ on the set of walks on $(Q,I)$, as the
smallest equivalence relation satisfying the following conditions :

\begin{enumerate}
\item For each arrow $\alpha :x\rightarrow y$, we have $\alpha \alpha
^{-1}\sim e_{x}$ and $\alpha ^{-1}\alpha \sim e_{y}$

\item For each minimal relation $\sum_{i=1}^{m}\lambda _{i}w_{i}$, we have $%
w_{i}\sim w_{j}$ for all $i$, $j$ in $\{1,\ldots ,m\}$.

\item If $u$, $v$, $w$ and $w^{\prime }$ are walks, and $u\sim v$ then $%
wuw^{\prime }\sim $ $wvw^{\prime }$, whenever these compositions are defined.
\end{enumerate}

We denote by $\tilde{w}$ the homotopic class of a walk $w$. Let $v_{0}$ be a
fixed point in $Q_{0}$, and consider the set $W(Q,v_{0})$ of walks of source
and target $v_{0}$. On this set, the product of walks is everywhere defined.
Because of the first and the third conditions in the definition of the
relation , one can form the quotient group $W(Q,v_{0})/\sim $ . This group
is called the fundamental group of the bound quiver $(Q,I)$ with base point $%
v_{0}$, denoted by $\pi _{1}(Q,I,v_{0})$. It follows from the connectedness
of $Q$ that this group does not depend on the base point $v_{0}$, and we
denote it simply by $\pi _{1}(Q,I)$. This group has a clear geometrical
interpretation as the first homotopic group of a C.W. complex $B(Q,I)$
associated to $(Q,I)$, see \cite{bustamante}.

\begin{remark}
If $I_{0}$ is null or monomial It is know that $\pi _{1}(Q,I_{0})$ is the
free product of $\chi (Q)$ copies of $%
\mathbb{Z}
$, where $\chi (Q)=\left\vert Q_{1}\right\vert $ $-\left\vert
Q_{0}\right\vert +1$ is the Euler characteristic of the underlying graph.
Thus, $\pi _{1}(Q,I_{0})= \mathbb{Z}^{\sqcup \chi (Q)}$.
\end{remark}

\begin{example}
\label{ejemplo1}Consider the quiver $O$ 
\begin{equation*}
\begin{array}{ccccc}
& u &  &  &  \\ 
b\nearrow &  & \searrow c &  &  \\ 
v & \underset{a}{\longrightarrow } & y & \underset{d}{\longrightarrow } & z%
\end{array}%
\end{equation*}%
and the clearly admissible ideal $I_{1}=\left\langle ad\right\rangle $. We
defined the morphism $\nu :kO\rightarrow kO/I_{1}\simeq A$ such that $%
v(a)=a-bc+I_{1}$,and $v(\gamma )=\gamma +I_{1}$ for any other path $\gamma $%
. We defined $I_{2}=\ker (\nu )$. By the Remark \ref{remark Tranversiones}, $%
\nu $ is a presentation of $A$ and therefore $kO/I_{1}\simeq kO/I_{2}\simeq
A $. On the other hand, since $I_{1}$ is generated by monomials relations
and $\chi (Q)=4-4+1$, we have $\pi _{1}(O,I_{1})=%
\mathbb{Z}
^{\sqcup 1}=%
\mathbb{Z}
$. Since $I_{2}=\left\langle ad+bcd\right\rangle $, we have $ad\sim bcd$ and 
$a\sim bc$. Therefore, it is clear that $\pi _{1}(O,I_{2})=\{1\}$.
\end{example}

\subsection{Co-products}

Given two pointed bounder quivers $(Q^{\prime },I^{\prime },v^{\prime })$
and $(Q^{\prime \prime },I^{\prime \prime },v^{\prime \prime })$ without any
vertex or arrow in common, we defined the bounded quiver $(Q,I):=(Q^{\prime
},I^{\prime },v^{\prime })\sqcup (Q^{\prime \prime },I^{\prime \prime
},v^{\prime \prime })$ in the following way: $Q_{0}:=Q_{0}^{\prime }\cup
Q_{0}^{\prime \prime }$ but we identify $v^{\prime }$ and $v^{\prime \prime
} $ to single vertex $v$, and $Q_{1}:=Q_{1}^{\prime }\cup Q_{1}^{\prime
\prime }$. Thus $Q^{\prime }$ and $Q^{\prime \prime }$ are identify with a
full convex sub-quivers of $Q$, and therefore $I^{\prime }$ and $I^{\prime
\prime }$ are ideals of $Q$. Finally we defined $I:=I^{\prime }+I^{\prime
\prime }$. Note that if $I^{\prime }$ and $I^{\prime \prime }$ are
admissible then $I$ as well. The vertex $v^{\prime }$ and $v^{\prime \prime
} $ are not relevant in general, therefore we will omit him and the word 
\textit{pointed}.

We will use the symbol $\sqcup $ as much to mention the co-product between
bounder quiver as mention the free product between groups. Note that this
operations are associative and commutative.

In \cite{Castonguay} is proved two important result for us. The first say
that the co-product behaves well under changes of presentations:

\begin{proposition}
\label{cambiobien}Let $A$ be an algebra with two presentations $%
(Q_{A},I_{A}) $ and $(Q_{A},I_{A}^{\prime })$ and let $B$ be an algebra with
two presentations $(Q_{B},I_{B})$ and $(Q_{B},I_{B}^{\prime })$. If 
\begin{eqnarray*}
(Q,I) &=&(Q_{A},I_{A})\sqcup (Q_{B},I_{B})\text{ and} \\
(Q,I^{\prime }) &=&(Q_{A},I_{A}^{\prime })\sqcup (Q_{B},I_{B}^{\prime })
\end{eqnarray*}%
then $(Q,I)$ and $(Q,I^{\prime })$ are presentation of the same algebra.
\end{proposition}

The second is analogous to the result to Van Kampen's theorem for
topological spaces:

\begin{proposition}
\label{grupobien}If $(Q,I)=(Q^{\prime },I^{\prime })\sqcup (Q^{\prime \prime
},I^{\prime \prime })$ then 
\begin{equation*}
\pi _{1}(Q,I)=\pi _{1}(Q^{\prime },I^{\prime })\sqcup \pi _{1}(Q^{\prime
\prime },I^{\prime \prime })
\end{equation*}
\end{proposition}

The two previous propositions and the next clear proposition will be very
used in our work.

\begin{proposition}
\label{triangularbien}If $(Q,I)=(Q^{\prime },I^{\prime })\sqcup (Q^{\prime
\prime },I^{\prime \prime })$ and the algebras $kQ^{\prime }/I^{\prime }$
and $kQ^{\prime \prime }/I^{\prime \prime }$ are triangular, then $kQ/I$ is
triangular.
\end{proposition}

\begin{remark}
\label{remark muiltiples usos} It is important to note that the propositions %
\ref{grupobien} and \ref{triangularbien} are easily generalized for a finite
amount of bounder quivers. In proposition \ref{cambiobien} we have two pair
of presentation of the same algebra respectively. This result can be
generalizable for any finite amount of pairs in this sense: If we have
collections of pairs $\{((Q_{A_{i}},I_{A_{i}}),(Q_{A_{i}},I_{A_{i}}^{\prime
}))\}_{i=1}^{n}$ such that the elements of each par is a presentation of the
same algebra $A_{i}$, then any co-product of $n$ presentations taking one,
and only one, presentation of each pair, is a presentation of the same
algebra.
\end{remark}

\begin{example}
\label{ejemploZ2}Consider the bounded quivers $(O^{\prime },I_{1}^{\prime })$
and $(O^{\prime \prime },I_{1}^{\prime \prime })$ isomorphous to the first
bounded quiver in the example \ref{ejemplo1}.The co-product $(O,I)$ between $%
(O^{\prime },I_{1}^{\prime },z^{\prime })$ and $(O^{\prime \prime
},I_{1}^{\prime \prime },v^{\prime \prime })$ is%
\begin{equation*}
\begin{array}{ccccccccc}
& u^{\prime } &  &  &  & u^{\prime \prime } &  &  &  \\ 
b^{\prime }\nearrow &  & \searrow c^{\prime } &  & b^{\prime \prime }\nearrow
&  & \searrow c^{\prime \prime } &  &  \\ 
v^{\prime } & \underset{a^{\prime }}{\longrightarrow } & y^{\prime } & 
\underset{d^{\prime }}{\longrightarrow } & t & \underset{a^{\prime \prime }}{%
\longrightarrow } & y^{\prime \prime } & \underset{d^{\prime \prime }}{%
\longrightarrow } & z^{\prime \prime }%
\end{array}%
\end{equation*}%
Note that by the Proposition \ref{grupobien}, $\pi _{1}(O,L)=%
\mathbb{Z}
\sqcup 
\mathbb{Z}
=%
\mathbb{Z}
^{\sqcup 2}$. On the other hand we make the co-product $(O,\bar{L})$ between 
$(O^{\prime },I_{2}^{\prime },z^{\prime })$ and $(O^{\prime \prime
},I_{2}^{\prime \prime },v^{\prime \prime })$ isomorphous to the second
bounded quiver in example \ref{ejemplo1}. Again, $\pi (O,\bar{L}%
)=\{1\}\sqcup \{1\}=\{1\}$. Like in example \ref{ejemplo1} $(O^{\prime
},I_{1}^{\prime })$ and $(O^{\prime },I_{2}^{\prime })$ are presentations of
the same algebra and $(O^{\prime \prime },I_{1}^{\prime \prime })$ and $%
(O^{\prime \prime },I_{2}^{\prime \prime })$ as well; hence we have the
hypothesis of proposition\ \ref{cambiobien}. Thus $(O,L)$ and $(O,\bar{L})$
are presentations of the same algebra. More, $O^{\prime }$ and $O^{\prime
\prime }$ are triangular, therefore, by proposition \ref{triangularbien}, $O$
is triangular. Summarizing: For the group $%
\mathbb{Z}
^{\sqcup 2}$ we got a triangular algebra $A:=kO/L$ with two presentation $%
(O,L)$ and $(O,\bar{L})$, such that $\pi _{1}(O,L)$ is $%
\mathbb{Z}
^{\sqcup 2}$ and $\pi _{1}(O,\bar{L})$ is trivial.
\end{example}

Is clearly that for get a bounded quiver with group $%
\mathbb{Z}
^{\sqcup m}$ and triangular algebra we can built by induction an analogous
construction of example \ref{ejemploZ2}. We can use more triangles and the
remark \ref{remark muiltiples usos}. Therefore we are in conditions to
statement the following lemma:

\begin{lemma}
\label{lemmaZm}For each $m$ natural, there exits a triangular algebra $A$
with two presentation $(O,L)$ and $(O,\bar{L})$ such that $\pi _{1}(O,L)=%
\mathbb{Z}
^{\sqcup m}$ and $\pi _{1}(O,\bar{L})=\{1\}.$
\end{lemma}

In the next section we will put our effort to prove the analogous for any\
finitely presented group instead of $%
\mathbb{Z}
^{\sqcup m}$.

\section{Fundamental group without free generator.\label{sec fgwfg}}

If $G$ is a finitely presented group, a \textbf{presentations} of $G$
consists of a finite set of generators $S$, and a finite set of relations
among these generators $R$ such that $G=\left\langle S~|R\right\rangle $
(see \cite{Hungerford}). If $x\in S$ but neither $x$ nor a power of $x$
appear in any relation of $R$, we say that $x$ is a \textbf{free generator}.
Let $\mathbb{H}$ be the set of all finitely presented group without free
generators. It is possible to write $G$ as free product between a group
generated by its free generator and a group in $\mathbb{H}$. In fact, we
consider $S_{F}=\{x\in S:x$ is a free generator$\}$, $m$ the amount of
element of $S_{F}$ and $H=\left\langle S~-S_{F}|R\right\rangle $. We have $G=%
\mathbb{Z}
^{\sqcup m}\sqcup H$. In this way, we divided the problem to find a algebra
with presentation whose group is $G$ in two problem: to find a algebra with
a presentation whose group is $%
\mathbb{Z}
^{\sqcup m}$ (solved by lemma \ref{lemmaZm}) and to find a algebra with a
presentation whose group is $H$. Hereafter we can handle to prove the next
key proposition:

\begin{proposition}
\label{prop principal}Let $H$ be a finitely presented without free
generators group, there exits a triangular algebra $A$ with two presentation 
$(Q_{H},I)$ and $(Q_{H},\bar{I})$ such that $\pi _{1}(Q_{H},I)=H$ and $\pi
_{1}(Q_{H},\bar{I})=\{1\}$.
\end{proposition}

The idea for the proof of \ref{prop principal} is the following. First, we
obtain $H$ generated by $w_{1},w_{2},\ldots ,w_{n}$ with two type of
relations: $w_{x}\cdots w_{x+t}=1$ (lasso) and $w_{a}=w_{b}$ (cross-lasso).
The quiver will be a chain of consecutive triangles, one for each generator $%
w_{i}$. The relations in the groups will be associated to the minimal
relations in $kQ_{H}$, that will generate $I$. Namely, a lasso relation $%
w_{x}\cdots w_{x+t}=1$ can be thought as a lasso which catches the triangles 
$x,x+1,\ldots ,x+t$. The cross-lasso relations $w_{a}=w_{b}$ can be thought
as a twist-lasso which catches the triangles $a$ and $b$. The idea is more
flexible than the one showed throughout this work, in fact, with this idea
it is possible to write a presentation with smaller (less vertices and
arrows) quivers that we present here.

\subsection{A particular presentation of a group\label{sub seca presentacion}%
}

Let $H$ be group in $\mathbb{H}$ with $H=\left\langle S~|R\right\rangle $.
For $H$ we will to give two special presentations $\left\langle S^{\prime
}~|R^{\prime }\right\rangle $ and $\left\langle S^{\prime \prime
}~|R^{\prime \prime }\right\rangle $ with special proprieties. The first
will be only auxiliary to get the second.

The general procedure to obtain the first alternative presentations is the
next. We consider an original presentation $\left\langle S~|R\right\rangle $%
. Without lost of generality we can suppose that each relations of $R$ is a
products of generators in the left side without powers and the neutral
element in the right side:%
\begin{eqnarray}
y_{1}^{1}\cdots y_{m_{1}}^{1} &=&1  \label{for ecuaciones primears} \\
&&\vdots  \notag \\
y_{1}^{s}\cdots y_{m_{s}}^{s} &=&1  \notag
\end{eqnarray}%
Note that two or more variables can represent the same generator. For
example, the group $\mathbb{Z}_{2}\oplus \ \mathbb{Z}$ has the presentations 
$\left\langle \{a,b\}~|\{aba^{-1}b^{-1}=1,~a^{2}=1\}\right\rangle $. We can
choice 
\begin{equation*}
\left\langle \{a,b,c,d\}~|\{abcd=1,ac=1,~bd=1,aa=1\}\right\rangle .
\end{equation*}%
Where $y_{3}^{1}=y_{2}^{2}=c$. We also impose that there are not relation
with only one variable, e.i., $y_{b}^{a}=1$.

Note that since $H\in \mathbb{H}$, \ every generator appear at least in one
relation of $R$.

The next step is created $S^{\prime }$ and $R^{\prime }$. The generators in $%
S^{\prime }$ will be $g_{1}^{1},\ldots ,\ldots ,g_{m_{s}}^{s}$, all
different. To construct the new relations in $R^{\prime }$ we proceed in
this way. In (\ref{for ecuaciones primears}) we replace $y_{k_{i}}^{i}$ by $%
g_{k_{i}}^{i}$ and, thus, we have another $s$ relations. 
\begin{eqnarray*}
g_{1}^{1}\cdots g_{m_{1}}^{1} &=&1 \\
&&\vdots \\
g_{1}^{s}\cdots g_{m_{s}}^{s} &=&1
\end{eqnarray*}

If some $y_{a}^{b}$ is the same element that $y_{c}^{d}$ , then we add the
relation $g_{a}^{b}=g_{c}^{d}$. Some of these relations can be redundant and
we can suppress these. Is evident that this presentation is equivalent to
the forward.

\begin{example}
\label{ejemplo ZZ2}For $H=\mathbb{Z}\oplus \mathbb{Z}_{2}$ we choice the
presentation%
\begin{equation*}
\left\langle \{a,b,c,d\}~|\{abcd=1,~ac=1,~bd=1,aa=1\}\right\rangle .
\end{equation*}%
The four relations produce the four new relations%
\begin{equation*}
g_{1}^{1}g_{2}^{1}g_{3}^{1}g_{4}^{1}=1,~g_{1}^{2}g_{2}^{2}=1,~g_{1}^{3}g_{2}^{3}=1,~g_{1}^{4}g_{2}^{4}=1.
\end{equation*}%
The aggregate relations are $g_{1}^{1}=g_{1}^{2}$,$g_{1}^{2}=g_{1}^{4}$,$%
g_{1}^{4}=g_{2}^{4}$ (for $a$),$g_{2}^{1}=g_{1}^{3}$(for $b$), $%
g_{3}^{1}=g_{2}^{2}$ (for $c$), $g_{4}^{1}=g_{2}^{3}$ (for $d$). Note that
the equality $g_{1}^{4}=g_{2}^{4}$ involved two elements of the same
relation. This fact is not convenient for the construction which we want to
do. Hence, we need a one more type of presentation.
\end{example}

For built $\left\langle S^{\prime \prime }|R^{\prime \prime }\right\rangle $
we procedure in this way: Let $n$ be the amount os variables in $S^{\prime }$%
. We simply rename all variables $g_{1}^{1},\ldots ,\ldots ,g_{m_{s}}^{s}$
by $w_{1},\ldots ,w_{n}$ respecting the order. Note that $n=m_{1}+\cdots
m_{s}$.This presentation is equivalent to the second (and, therefore,
equivalent to the original) because we only made a rename of variables.

The specials proprieties which $R^{\prime \prime }$ has are:

\begin{description}
\item[(a)] Note that all the relations in $R^{\prime \prime }$ has the form $%
w_{i_{1}}\cdots w_{i_{k}}=1$ (we will call \textit{lasso}) or $w_{i}=w_{j}$ (%
\textit{cross-lasso}).

\item[(b)] Each $w_{i}$ appear at most ones in each lasso and cross-lasso
relations.

\item[(c)] In the relations $w_{i_{1}}\cdots w_{i_{k}}=1$ the index are
consecutive, e.i., $i_{k+1}=i_{k}+1$, therefore this relations has the form $%
w_{i},w_{i+1}\ldots ,w_{i+k}=1$
\end{description}

\section{The quivers and the ideals\label{sec quiver and ideal}}

Let $H\in \mathbb{H}$. We consider the presentation $\left\langle S^{\prime
\prime }~|R^{\prime \prime }\right\rangle $ made in the section \ref{la
presentacion}. Based on this presentation we will build a quiver $Q$ and a
ideal $I$ such that the homotopic group of $(Q_{H},I)$ will be $H$. The
quiver is based in $S^{\prime \prime }$ and the ideal on $R^{\prime \prime }$%
. After we will build another ideal $\bar{I}$ such that $\pi _{1}(Q_{H},\bar{%
I})=\{1\}$.

\subsection{The quiver\label{subsec The quiver}}

We defined the quiver $Q_{H}:=(Q_{0},Q_{1})$ where 
\begin{equation*}
Q_{0}=\{x_{1},\ldots ,x_{n+1},y_{1},\ldots ,y_{n}\}
\end{equation*}%
and $Q_{1}$ is formed by following collections of arrows%
\begin{eqnarray*}
r_{i} &:&x_{i}\rightarrow y_{i}\text{ for }i=1,\ldots ,n \\
l_{i} &:&y_{i}\rightarrow x_{i+1}\text{ for }i=1,\ldots ,n \\
a_{i} &:&x\rightarrow x_{i+1}\text{ for }i=1,\ldots ,n.
\end{eqnarray*}%
Remark: Let $n$ be the amount of elements of $S^{\prime \prime }$. The
quiver $Q_{H}$ is%
\begin{equation*}
\begin{array}{ccccccccccc}
& y_{1} &  &  &  & y_{i} &  &  &  & y_{n} &  \\ 
r_{1}\nearrow &  & \searrow l_{1} & \cdots & r_{i}\nearrow &  & \searrow
l_{1} & \cdots & r_{n}\nearrow &  & \searrow l_{n} \\ 
x_{1} & \underset{a_{1}}{\longrightarrow } & x_{2} &  & x_{i} & \underset{%
a_{i}}{\longrightarrow } & x_{i+1} &  & x_{n} & \underset{a_{n}}{%
\longrightarrow } & x_{n+1}%
\end{array}%
\end{equation*}%
We can described this as a chain of $n$ triangles.

\begin{notation}
We take the following notations%
\begin{eqnarray*}
m_{ij} &:&=a_{i+1}\cdots a_{j-1} \\
A_{i} &:&=r_{i}l_{i} \\
T_{i} &:&=A_{i}a_{i}^{-1} \\
\alpha _{i} &:&=m_{0i}T_{i}m_{0i}^{-1}
\end{eqnarray*}
\end{notation}

\subsection{The ideals and the homotopy \label{sebsec the ideals}}

\subsubsection{The homotopy of the bounder quivers\label{subsubsec The homot
of}}

We will study the homotopy of the quiver $Q_{H}$ defined in \ref{el quiver}.
Consider the null ideal $I_{0}$. Note that $\chi (Q_{H})=$ $3n-(2n+1)-1=n$.
We choice the base point $x_{1}$. Thus $\pi _{1}(Q_{H},I_{0},x_{1})=%
\mathbb{Z}
^{\sqcup n}$. Note that the walks $\alpha _{i}$ are loops with base in $%
x_{1} $. Applying technical of algebraic topology to the quiver, it is clear
that each class $\tilde{\alpha}_{i}$ is generator of $\pi
_{1}(Q_{H},I_{0},x_{1})$, i.e., $\pi _{1}(Q_{H},I_{0},x_{1})=$ $\left\langle
\{\tilde{\alpha}_{1},\ldots ,\tilde{\alpha}_{n}~|\varnothing \right\rangle $%
. More, if we now consider a new ideal $I$ in $kQ_{H}$, then $\pi
_{1}(Q_{H},I_{0},x_{1})=\left\langle \{\tilde{\alpha}_{1},\ldots ,\tilde{%
\alpha}_{n}~|R\right\rangle $ where the relations in $R$ must be generated
by homotopic relations generates by the minimal relations in $I$. We will
rely on this fact to make the computation of the fundamental groups of the
bounder quivers of the subsection \ref{subsubsec the ideal for h}.

\subsubsection{The ideal for $H$.\label{subsubsec the ideal for h}}

We will express the ideal $I$ generated by relations in the algebra $kQ_{H}$%
. For each relations in $R^{\prime \prime }$ we consider a relations in $I$
in the next way. If the relations in $R^{\prime \prime }$ has the form $%
w_{i}w_{i+1}\cdots w_{j}=1,$ then we consider the relation $%
a_{i}a_{i+1}\cdots $ $a_{j}+A_{i}A_{i+1}\cdots $ $A_{j}\ $(we call \textit{%
lasso} because the induced homotopy is $A_{i}A_{i+1}\cdots $ $%
A_{j}a_{j}^{-1} $ $\cdots a_{i+1}^{-1}a_{i}\sim 1$, which can be thought as
a lasso that envelops the triangles $T_{i},\ldots T_{j}$). Note that this it
is possible because the index of $w_{i}w_{i+1}\cdots w_{j}$ are consecutive
(property (c)). On a other hand, if the relations has the form $w_{i}=w_{j}$
with $i<j$ then we put $a_{i}m_{ij}A_{j}$ $+A_{i}m_{ij}a_{j}$ (we call 
\textit{cross-lasso }because the induced homotopy is $a_{i}m_{ij}A_{j}$ $%
\sim A_{i}m_{ij}a_{j}$ , which can be thought as a lasso that envelops the
triangles $T_{i}$ and $T_{j}$).

Note that we never has a monomial relation as $a_{i}+A_{i}$, therefore we
have $I\subset F^{2}$. And, since $Q_{H}$ has not oriented cycled for some $%
m $, one has $F^{m}=0\subset I$. Therefore $I$ is admissible.

\begin{example}
Continuing with $\mathbb{Z}\oplus \mathbb{Z}_{2}$. The elements of $%
S^{\prime \prime }$ are $w_{1},\ldots ,w_{10}$. Hence the quiver has ten
triangles.$~$For example, the relations $w_{1}w_{2}w_{3}w_{4}=1$ in $%
R^{\prime \prime }$ is in correspondence with the relation $%
a_{1}a_{2}a_{3}a_{4}+A_{1}A_{2}A_{3}A_{4}$ in $kQ_{H}$, and the relation $%
w_{1}=w_{5}$ corresponding to $a_{1}m_{1,5}A_{5}$ $+A_{1}m_{1,5}a_{5}$. The
ideal $I$ is generated by this relations and the other six of example \ref%
{ejemplo ZZ2}.
\end{example}

Now we will to prove that all this relations together are minimal. The
unique sub-algebras which has some relations, are $I(x_{i}x_{j})$. And the
unique possible relation in this sub-algebra are $a_{i}m_{ij}A_{j}$ $%
+A_{i}m_{ij}a_{j}$ or $a_{i}\cdots $ $a_{j}+A_{i}\cdots $ $A_{j}$. Note that
the four terms involved are different, therefore, the two relations are
minimal in the ideals generated by themselves.

It remains to prove that the homotopic group of $(Q_{H},I)$ is isomorphous
to $H$. For this we will to prove that the classes of the loops $\tilde{%
\alpha}_{1},\ldots ,\tilde{\alpha}_{n}$ check exactly the same relations
that the elements $w_{1}\ldots ,w_{n}$ check in $H$.

\begin{lemma}
The relation $a_{i}m_{ij}A_{j}$ $+A_{i}m_{ij}a_{j}$ generates in $\pi
_{1}(Q_{H},I,x_{1})$ the homotopic relation $\tilde{\alpha}_{i}\sim \tilde{%
\alpha}_{j}$.
\end{lemma}

\begin{proof}
Since $a_{i}m_{ij}A_{j}$ $+A_{i}m_{ij}a_{j}$ is minimal in $I$ then $%
a_{i}m_{ij}A_{j}\sim A_{i}m_{ij}a_{j}$. For simplicity we denote $m_{ij}$ by 
$m$. Therefore%
\begin{eqnarray*}
A_{i}ma_{j} &\sim &a_{i}mA_{j} \\
A_{i}a_{i}^{-1}a_{i}ma_{j} &\sim &a_{i}mA_{j} \\
A_{i}a_{i}^{-1} &\sim &a_{i}mA_{j}a_{j}^{-1}m^{-1}a_{i}^{-1} \\
T_{i} &\sim &a_{i}mT_{j}(ma_{i})^{-1} \\
m_{0i}T_{i}m_{0i}^{-1} &\sim &m_{0i}a_{i}mT_{j}(ma_{i})^{-1}m_{0i}^{-1} \\
\alpha _{i} &\sim &\alpha _{j}
\end{eqnarray*}
\end{proof}

\begin{lemma}
The relation $a_{i}a_{i+1}\cdots $ $a_{i+k}+A_{i}A_{i+1}\cdots $ $A_{i+k}$
generates in $\pi _{1}(Q_{H},I,x_{1})$ the homotopic relation $\tilde{\alpha}%
_{i}\tilde{\alpha}_{i+1}\cdots \tilde{\alpha}_{i+k}\sim 1$.
\end{lemma}

\begin{proof}
For simplicity, we denote $a_{i+j}$, $A_{i+j}$, $m_{0i}$ by $b_{j}$, $B_{j}$
and $m$, respectively. Since $b_{0}\cdots $ $b_{k}+B_{0}\cdots $ $B_{k}$ is
minimal in $I$ then 
\begin{eqnarray*}
B_{0}\cdots B_{k} &\sim &b_{0}\cdots b_{k} \\
B_{0}\cdots B_{k}(b_{0}\cdots b_{k})^{-1} &\sim &e_{x_{i}} \\
mB_{0}\cdots B_{k}(b_{0}\cdots b_{k})^{-1}m^{-1} &\sim &e_{x_{1}} \\
mB_{0}(mb_{0})^{-1}(mb_{0})B_{1}\cdots B_{k}(mb_{1}\cdots b_{k})^{-1} &\sim
&e_{x_{1}}
\end{eqnarray*}%
Note that $mB_{0}(mb_{0})^{-1}\sim \alpha _{i}$. Thus,%
\begin{eqnarray*}
\alpha _{i}(mb_{0})B_{1}\cdots B_{k}(mb_{0}\cdots b_{k})^{-1} &\sim
&e_{x_{1}} \\
\alpha _{i}(mb_{0})B_{1}(mb_{0}b_{1})^{-1}(mb_{0}b_{1})\cdots
B_{k}(mb_{0}\cdots b_{k})^{-1} &\sim &e_{x_{1}}
\end{eqnarray*}%
Note that $(mb_{0})B_{1}(mb_{0}b_{1})^{-1}\sim \alpha _{i+1}$. Hence we
repeat the procedure until the last step%
\begin{eqnarray*}
\alpha _{i}\alpha _{i+1}\cdots \alpha _{i+k-1}(mb_{0}b_{1}\cdots
b_{k-1})B_{k}(mb_{0}\cdots b_{k})^{-1} &\sim &e_{x_{1}} \\
\alpha _{i}\alpha _{i+1}\cdots \alpha _{i+k-1}\alpha _{i+1} &\sim &e_{x_{1}}
\end{eqnarray*}
\end{proof}

Finally, note that for each relation in $R^{\prime \prime }$ of the kind $%
w_{i}\cdots w_{j}=1$, we have $\tilde{\alpha}_{i}\cdots \tilde{\alpha}%
_{j}\sim 1$ in $\pi _{1}(Q_{H},I,x_{1})$. And for each relation in $%
R^{\prime \prime }$ of the kind $w_{i}=w_{j}$, we have $\tilde{\alpha}_{i}=%
\tilde{\alpha}_{j}$ in $\pi _{1}(Q_{H},I,x_{1})$. Therefore $H=\pi
_{1}(Q_{H},I,x_{1})=\pi _{1}(Q_{H},I)$. Thus we proofed a first part of the
Proposition \ref{prop principal}.

\subsubsection{The ideal for the trivial group.\label{subsubsec the ideal
for 1}}

The new ideal is the image of the previous ideal $I$ by composition of
trasversions. Actually, each transvection $\varphi _{j}$ \textit{kill} the
element $\tilde{\alpha}_{i}$ in the fundamental group $\pi
_{1}(Q_{H},\varphi _{j}(I))$. We not will to prove exactly that, however for
alternative understanding and realization of this, it can be interested see 
\cite{patric}.

\begin{remark}
If $\bar{\eta}:kQ_{H}\rightarrow kQ_{H}$ is an automorphism of algebra such
that $\bar{\eta}^{-1}$ is a transversion and $J$ is a admissible ideal, then 
$(Q_{H},\eta (J))$ and $(Q_{H},J)$ are presentations of the same algebra,
because, when consider the projection \ $\eta :kQ_{H}\rightarrow kQ_{H}/J$,
we have that $\bar{\eta}(J)=\ker \eta ^{-1}$.
\end{remark}

\begin{notation}
Let $s$ be the amount of lassos. If the lasso $a_{i}\cdots $ $%
a_{i+h}+A_{i}\cdots $ $A_{i+h}$ is the $k$-th lasso, we rename the walks
such that this relation result $a_{1}^{k}\cdots $ $a_{n_{k}}^{k}+A_{1}^{k}%
\cdots $ $A_{n_{k}}^{k}$, we abbreviate this by $W^{k}$. Let $\varphi
_{_{i}}^{k}:kQ_{H}\rightarrow kQ_{H}$ be automorphism of algebras defined by 
$\varphi _{_{i}}^{k}(a_{_{i}}^{k})=a_{_{i}}^{k}+A_{_{i}}^{k}$ and the
identity at other vertex and arrow. We defined $\bar{\gamma}^{k}:=\varphi
_{n_{k}}^{k}\circ \cdots \circ \varphi _{1}^{k}$ and $\bar{\gamma}=$ $\bar{%
\gamma}^{s}\circ \cdots \circ \bar{\gamma}^{1}$.
\end{notation}

\begin{lemma}
\label{lemmaTecnico}With the above notation%
\begin{equation}
\bar{\gamma}(W^{k})=a_{1}^{k}\cdots
a_{n_{k}}^{k}+\sum_{i=1}^{n_{k}}a_{1}^{k}\cdots
a_{i-1}^{k}A_{i}^{k}a_{i+i}^{k}\cdots a_{n_{k}}^{k}+\Gamma _{k}
\label{for kill relation}
\end{equation}%
where $\Gamma _{k}$ is a sum such that each term has two or more capital
letters.
\end{lemma}

\begin{proof}
Note $\bar{\gamma}(W^{k})=$ $\bar{\gamma}^{s}\circ \cdots \circ \bar{\gamma}%
^{1}(W^{k})=\bar{\gamma}^{k}(W^{k})$, because $\bar{\gamma}_{k}$ the
identity over $I(x_{1}^{k},x_{n_{k}-1}^{k})$ when $l\neq k$. We fix a $k$
and for to make more clear the notations we replace: $n_{k}$ for $n$, $a$ by 
$b$ and $A$ by $B$, and we omit the upper-index $k$ in the letters $b$, $B$
and $\varphi $. Thus the inductive hypothesis over $i$ is%
\begin{equation}
\varphi _{i}\circ \cdots \circ \varphi _{1}(W^{k})=b_{1}\cdots
b_{n}+\sum_{j=1}^{i}b_{1}\cdots B_{j}\cdots b_{n}+\Gamma _{k}^{i}
\label{sumaFea}
\end{equation}%
where $\Gamma _{k}^{i}$ is a sum such that each term has two or more capital
letters. For $i=1$ is clearly, in fact $\varphi _{1}(W^{k})=b_{1}\cdots
b_{n}+B_{1}b_{2}\cdots b_{n}+B_{1}\cdots B_{n}$, and $n\geq 2$. Suppose that
the hypothesis is true for $i$. Applying $\varphi _{i+1}$ in (\ref{sumaFea})
and reorganizing we have: 
\begin{eqnarray*}
&=&b_{1}\cdots b_{n}+b_{1}\cdots B_{i+1}\cdots
b_{n}+\sum_{j=1}^{i}(b_{1}\cdots B_{j}\cdots b_{n}+b_{1}\cdots B_{j}\cdots
B_{i+1}\cdots b_{n}) \\
&&+\varphi _{i+1}(\Gamma _{k}^{i}) \\
&=&b_{1}\cdots b_{n}+\sum_{j=1}^{i}(b_{1}\cdots B_{j}\cdots
b_{n})+b_{1}\cdots B_{i+1}\cdots b_{n}+ \\
&&+\sum_{j=1}^{i}(b_{1}\cdots B_{j}\cdots B_{i+1}\cdots b_{n})+\varphi
_{i+1}(\Gamma _{k}^{i}) \\
&=&b_{1}\cdots b_{n}+\sum_{j=1}^{i+1}(b_{1}\cdots B_{j}\cdots b_{n})+\Gamma
_{k}^{i+1}
\end{eqnarray*}%
where it is clear that $\Gamma _{k}^{i+1}$is a sum such that each term has
two or more capital letters because: $\Gamma _{k}^{i}$ has the same
propriety and when $\varphi _{i+1}$ is applied over it, the amount of
capital letter can never decrease, and, on another hand, it is clearly that $%
\Sigma _{i=1}^{j}(b_{1}\cdots B_{j}\cdots B_{i+1}\cdots b_{n})$ has two
capital letter. So, applying induction until $n$, we have that $\bar{\gamma}%
^{k}(W^{k})=\varphi _{n}\circ \cdots \circ \varphi _{1}(W^{k})$ has the
required form.
\end{proof}

Now we can finish the proof of the proposition \ref{prop principal}.

\begin{proof}
Defined $\bar{I}:=\bar{\gamma}(I)$. It is easy check that $\bar{I}$ is
admissible as well. Since $\bar{\gamma}$ is automorphism, the minimal
relation in $I$ are minimal relations in $\bar{I}$, in particular the
relation (\ref{for kill relation}). By the lemma we can ensure that the
homotopy%
\begin{equation*}
a_{1}^{k}\cdots a_{i}^{k}\cdots a_{n_{k}}^{k}\sim a_{1}^{k}\cdots
A_{i}^{k}\cdots a_{n_{k}}^{k}
\end{equation*}%
is verified. And from this we have $a_{i}^{k}\sim A_{i}^{k}$ for all $%
k=1,\ldots ,s$ and $i=1,\ldots ,n_{k}$. Since all $w_{j}$ is in some lasso,
we have $a_{j}\sim A_{j}$ for all $j=1,\ldots n$. Manipulating this, we get $%
e_{x_{1}}\sim m_{0j}T_{j}m_{0j}^{-1}\sim \alpha _{j}$. Thus, the fundamental
group of $\pi _{1}(Q_{H},\bar{I},x_{1})$ verify $\tilde{\alpha}_{j}\sim 1$
for all $j$. But $\{\tilde{\alpha}_{1},\ldots \tilde{\alpha}_{n}\}$ is a set
of generators, therefore $\pi _{1}(Q_{H},\bar{I})\approx \{1\}$.
\end{proof}

\subsection{The presentation for $G$ \label{subsec presentation for G}}

Let $G$ a finitely presented group such that $G=%
\mathbb{Z}
^{\sqcup m}\sqcup H$ with $H\in \mathbb{H}$. We defined the bounder quivers 
\begin{eqnarray*}
(Q_{G},J_{G}) &:&=(O,L)\sqcup (Q,I) \\
(Q_{G},\bar{J}) &:&=(O,\bar{L})\sqcup (Q,\bar{I})
\end{eqnarray*}%
where $O$, $L$ and $\bar{L}$ are the same of Lemma \ref{lemmaZm} and $Q$, $I$
and $\bar{I}$ are the same that in the Proposition \ref{prop principal}. By
Proposition \ref{cambiobien} these are presentation of the same algebra. By
Proposition \ref{grupobien} 
\begin{equation*}
\pi _{1}(Q_{G},J_{G})=\pi _{1}(O,L)\sqcup \pi _{1}(Q,I)=H\sqcup 
\mathbb{Z}
^{\sqcup m}=G
\end{equation*}%
and $\pi _{1}(Q_{G},\bar{J}_{G})=\{1\}\sqcup \{1\}=\{1\}$. By Proposition %
\ref{triangularbien} the quiver $Q_{G}$ is triangular. Thus we proved the
next proposition.

\begin{proposition}
Let $G$ be a finitely presented group, there exits a triangular algebra $A$
with two presentation $(Q_{G},J_{G})$ and $(Q_{G},\bar{J}_{G})$ such that $%
\pi _{1}(Q_{G},J_{G})=G$ and $\pi _{1}(Q_{G},\bar{J}_{G})=\{1\}$.
\end{proposition}

Now we prove the Theorem \ref{theo unico}.

\begin{proof}
Let $G_{1},\ldots ,G_{n}$ be a finitely presented groups. By the formula (%
\ref{for amalgamagrande}) we defined $(\hat{Q},I_{i})$. Similarly to do with 
$(Q_{G},J_{G})$ and $(Q_{G},\bar{J}_{G})$, by the remark \ref{remark
muiltiples usos}, we have, first:%
\begin{eqnarray*}
\pi _{1}(\hat{Q},I_{i}) &=&(\sqcup _{j\neq i}\pi _{1}(Q_{G_{j}},\bar{J}%
_{j}))\sqcup \pi _{1}(Q_{G_{i}},J_{i}) \\
&=&(\sqcup _{j\neq i}\{1\})\sqcup G_{i}=G_{i}
\end{eqnarray*}%
Second: we have pairs of presentation $((Q_{G_{i}},J_{G_{i}}),(Q_{G_{i}},%
\bar{J}_{G_{i}}))_{i=1}^{n}$ of the same algebras. Each $(\hat{Q},I_{i})$ is
made taking one\ and only one, presentation of each par, hence these are
presentations of the same algebra $A$. Third, all quivers involved in the
co-product are triangular, therefore the quivers $\hat{Q}$ is triangular,
and thus $A$ is triangular.
\end{proof}

\end{document}